\documentclass[12pt]{amsart}
\usepackage{mathrsfs}
\usepackage{amssymb}

\usepackage{titletoc}
\input xy
  \xyoption{all}
\usepackage{amscd}
\usepackage{amsmath, amssymb}
\usepackage{amsfonts}
\usepackage[colorlinks,linkcolor=blue,citecolor=blue, pdfstartview=FitH]{hyperref}

\numberwithin{equation}{section}

\theoremstyle{plain}
\newtheorem{thm}{Theorem}[section]
\newtheorem{theorem}[thm]{Theorem}
\newtheorem{lemma}[thm]{Lemma}
\newtheorem{corollary}[thm]{Corollary}
\newtheorem{proposition}[thm]{Proposition}

\theoremstyle{definition}

\newtheorem{question}[thm]{Question}
\newtheorem{problem}[thm]{Problem}
\newtheorem{remark}[thm]{Remark}

\newtheorem{definition}[thm]{Definition}

\newtheorem{example}[thm]{Example}

\newtheorem{defn-thm}[thm]{Definition-Theorem}

\newcommand{\sO}{{\mathcal O}}

\newcommand{\C}{{\mathbb C}}

\renewcommand{\P}{{\mathbb P}}
\newcommand{\Q}{{\mathbb Q}}
\newcommand{\R}{{\mathbb R}}

\newcommand{\bp}{\bar{\partial}}

\newcommand{\ds}{\oplus}

\newcommand{\ts}{\otimes}

\newcommand{\btheorem}{\begin{theorem}}
\newcommand{\etheorem}{\end{theorem}}
\newcommand{\bproposition}{\begin{proposition}}
\newcommand{\eproposition}{\end{proposition}}
\newcommand{\bdefinition}{\begin{definition}}
\newcommand{\edefinition}{\end{definition}}
\newcommand{\bcorollary}{\begin{corollary}}
\newcommand{\ecorollary}{\end{corollary}}
\newcommand{\bproof}{\begin{proof}}
\newcommand{\eproof}{\end{proof}}
\newcommand{\bremark}{\begin{remark}}
\newcommand{\eremark}{\end{remark}}
\newcommand{\eexample}{\end{example}}
\newcommand{\bexample}{\begin{example}}

\newcommand{\elemma}{\end{lemma}}
\newcommand{\blemma}{\begin{lemma}}

\newcommand{\sq}{\sqrt{-1}}

\newcommand{\p}{\partial}

\renewcommand{\bar}{\overline}
\newcommand{\eps}{\varepsilon}

\renewcommand{\phi}{\varphi}

\newcommand{\ee}{\end{eqnarray*}}
\newcommand{\be}{\begin{eqnarray*}}

\newcommand{\beq}{\begin{equation}}
\newcommand{\eeq}{\end{equation}}

\newcommand{\bd}{\begin{enumerate}}
\newcommand{\ed}{\end{enumerate}}

\renewcommand{\tilde}{\widetilde}

\newcommand{\rw}{\rightarrow}

\renewcommand{\rm}{\textrm}

\renewcommand{\>}{\rightarrow}

\usepackage{fancyhdr}
\renewcommand{\leftmark}{{\scriptsize RC-positivity and scalar-flat metrics on ruled surfaces}}

\begin{document}
\title{RC-positivity and scalar-flat metrics on ruled surfaces} \makeatletter
\let\uppercasenonmath\@gobble
\let\MakeUppercase\relax
\let\scshape\relax
\makeatother

\author{Jun Wang}
\date{}
\address{{
        Academy of Mathematics and Systems Science,
        Chinese Academy of Sciences, Beijing, 100190, China.}}

\author{Xiaokui Yang}
\date{}
\address{{Address of Xiaokui Yang: Morningside Center of Mathematics, Institute of
        Mathematics, Hua Loo-Keng Center of Mathematical Sciences,
        Academy of Mathematics and Systems Science,
        Chinese Academy of Sciences, Beijing, 100190, China.}}
\email{\href{mailto:xkyang@amss.ac.cn}{{xkyang@amss.ac.cn}}}

\noindent\thanks{This work was partially supported   by China's
Recruitment
 Program of Global Experts
 }
\maketitle

\begin{abstract}  Let $X$ be a ruled surface over a curve of genus $g$. We prove that   $X$ has a scalar-flat
Hermitian metric \emph{if and only if} $g\geq 2$ and $m(X)>2-2g$
where $m(X)$ is an intrinsic number depends  on the
complex structure of $X$. 
\end{abstract}
\setcounter{tocdepth}{1} \tableofcontents

\section{Introduction}

  In his ``Problem section", S.-T. Yau proposed the following
classical
 problem (\cite[Problem~41]{Yau82}), which is investigated
 intensively  in the last forty years.

\begin{problem}\label{problem} Classify all compact K\"ahler surfaces with zero scalar
curvature.
\end{problem}

\noindent By the celebrated Calabi-Yau Theorem (\cite{Yau78}), all
K\"ahler surfaces with vanishing first Chern class (e.g. $K3$
surfaces) admit K\"ahler metrics with zero scalar curvature.
 Such metrics are usually called \emph{scalar-flat}
K\"ahler metrics and it is a special class of constant scalar
curvature K\"ahler (cscK) metrics or extremal metrics. Obstructions
to the existence of such metrics have been known since the
pioneering works of S.-T. Yau \cite{Yau74} and E. Calabi
\cite{Cal85}. For comprehensive discussions on this rich topic, we
refer to \cite{Yau74, Yau78, Fut83, BD88, Tian90, Sim91, Fuj92,
LS93, LS94, Tian97, Don01, RS05, AP06, AT06, RT06, Ross06, CT08,
Sto08, AP09, ACGT11, Sze14, Sze17} and the references therein.

 In this paper, we study the geometry of compact complex manifolds with scalar-flat Hermitian
 metrics
 (with respect to the Chern connection), which is a generalization of
  Problem \ref{problem}. We begin with a characterization
 of compact complex manifolds with scalar-flat Hermitian metrics, which can be regarded as a
Hermitian analogue of Kazdan-Warner-Bourguignon's classical work in
Riemannian geometry, and we refer to
 \cite{Bes86} and \cite{Fut93} for more details.

\btheorem\label{main1} A compact complex manifold $X$ admits a
scalar-flat Hermitian metric if and only if $X$ is Chern Ricci-flat,
or both $K_X$ and $K_X^{-1}$ are RC-positive. \etheorem

\noindent Recall that, a line bundle $\mathscr L$ is called
\emph{RC-positive} if it has a smooth Hermitian metric $h$ such that
its curvature $-\sq\p\bp\log h$ has at least one positive eigenvalue
everywhere. By using  a remarkable theorem in \cite{TW10}
established by Tosatti-Weinkove (which is a Hermitian analogue of
Yau's theorem \cite{Yau78}), the anti-canonical bundle $K^{-1}_X$ is
RC-positive if and only if $X$ has a smooth Hermitian metric
$\omega$ such that its Ricci curvature $\mathrm{Ric}(\omega)$ has at
least one positive eigenvalue everywhere. A complex manifold $X$ is
called \emph{Chern Ricci-flat} if there exists a smooth Hermitian
metric $\omega$ such that the Chern-Ricci curvature
$\mathrm{Ric}(\omega)=-\sq\p\bp\log\omega^n=0$. On the other hand,
we proved in \cite[Theorem~1.4]{Yang17D} that a line bundle
$\mathscr L$ is RC-positive if and only if its dual line bundle
$\mathscr L^*$ is not pseudo-effective. By taking this advantage, we
can verify the RC-positivity of $K_X$ or $K_X^{-1}$ by adapting
methods in
differential geometry  as well as algebraic geometry. \\

As a straightforward application of Theorem \ref{main1}, we obtain

\bcorollary \label{main0} Let $X$ be a compact K\"ahler manifold. If
$X$ has a scalar-flat K\"ahler metric $\omega$, then either $X$ is a
Calabi-Yau manifold, or both $K_X$ and $K_X^{-1}$ are RC-positive.
\ecorollary

\noindent For instance, if $X$ is the blowing-up of $\P^2$ along
$m$-points ($m\leq 9$), it is well-known that the anti-canonical
bundle $K^{-1}_X$ is effective (e.g. \cite[p.~125-p.~129]{Fri98})
and so it is pseudo-effective. In this case, $K_X$ can not be
RC-positive and $X$ has no scalar-flat Hermitian (or K\"ahler)
metrics.

\bcorollary\label{P2} Let $\P^2\#m\bar{\P^2}$ be the blowing-up of
$\P^2$ along $m$ points. If $X$ admits a scalar-flat Hermitian
metric, then $m\geq 10$. \ecorollary

\noindent Indeed, it is proved by Rollin-Singer in
\cite[Theorem~1]{RS05} (see also \cite{Leb86, Leb91, LS93}) that: a
complex surface $X$ obtained by blowing-up $\P^2$ at $10$ suitably
chosen points admits a scalar-flat K\"ahler metric and any further
blowing-up of $X$ also admits a scalar-flat K\"aler metric.

A compact complex surface $X$ is called a \emph{ruled surface} if it
is a holomorphic $\P^1$-bundle over a compact Riemann surface $C$.
It is well known that any ruled surface $X$ can be written as a
projective bundle $\P(\mathscr E)$ where $\mathscr E$ is a rank two
vector bundle over $C$. Moreover, two ruled surfaces $\P(\mathscr
E)$ and $\P(\mathscr E')$ are isomorphic if and only if $\mathscr
E\cong \mathscr E'\otimes \mathscr L$ for some line bundle $\mathscr
L$ over $C$. The existence of cscK metrics on ruled surfaces are
extensively
studied, 
and we refer to \cite{Yau74, BD88, Tian90, Sim91, Fuj92, LS93, LS94,
RS05, AP06, AT06, RT06, Ross06, Sto08, ACGT11, Sze14} and the
references therein. A remarkable result (e.g. \cite{AT06, BD88,
ACGT11}) asserts that:
  \emph{A ruled surface $\P(\mathscr E)$ admits a cscK
metric if and only if $\mathscr E$ is poly-stable.}\\

 In the following, we aim to classify ruled surfaces with
scalar-flat Hermitian metrics. Let $\mathscr E$ be a rank two vector
bundle over a smooth curve $C$. One can define a number $m(\mathscr
E)$ (e.g. \cite[p.~122]{Fri98}) which is equal to  the \emph{minimal
degree} of $\mathscr E\otimes \mathscr L$ if there exists a
\textbf{sheaf extension} of $\mathscr E\otimes \mathscr L$: $$ 0\rw
\sO_C\rw \mathscr E\otimes \mathscr L\rw \mathscr F\rw0 $$ for some
line bundle $\mathscr L$. It is obvious that $m(\mathscr
E)=m(\mathscr E\ts \tilde {\mathscr L})$ for any line bundle $\tilde
{\mathscr L}$. Hence,  we can define an intrinsic number $m(X)$ for
a ruled surface $X$: $m(X)=m(\mathscr E)$ if $X$ can be written as
$\P(\mathscr E)$. It is obvious that $m(X)$ is independent of the
choices of $\mathscr E$. Let's explain the geometric meaning of
$m(X)$ by the example $X=\P(\mathscr L\ds \sO_C)\>C$ where $\mathscr
L$ is a line bundle. In this case, $m(X)=-|\deg(\mathscr L)|\leq 0$.
As another application of Theorem \ref{main1}, we obtain

\btheorem \label{main2} Let $X$ be a ruled surface over a smooth
curve $C$ of genus $g$. Then  $X$ has a scalar-flat Hermitian metric
if and only if $g\geq 2$ and $m(X)>2-2g$. \etheorem

\noindent  In particular, we have

\bcorollary\label{simple} Let $\mathscr L\>C$ be a line bundle over
a smooth curve of genus $g$ and $X=\P(\mathscr L\ds \sO_C)$. Then
$X$ has a scalar-flat Hermitian metric if and only if $g\geq 2$ and
$|\deg(\mathscr L)|< 2g-2$. \ecorollary

\noindent For instance, if $C$ is a smooth curve of degree $d>4$ in
$\P^2$, then the genus of $C$ is $g=\frac{1}{2}(d-1)(d-2)$ and the
degree of $\sO_C(1)$ is $d<2g-2$. Hence, $X=\P(\sO_C(1)\ds \sO_C)$
has scalar-flat Hermitian metrics. Note also that, in Corollary
\ref{simple}, if $\deg(\mathscr L)=0$, the vector bundle $\mathscr
L\ds \sO_C$ is poly-stable and $X=\P(\mathscr L\ds \sO_C)$ admits
scalar-flat K\"ahler metrics; however, when $0<|\deg(\mathscr
L)|<2g-2$, it  has no scalar-flat K\"ahler metrics. Moreover, we
construct such examples in higher dimensional ruled manifolds.

 \bproposition \label{main01}  Let $C$ be a smooth curve with genus $g\geq
 2$ and $\mathscr L$ be a line bundle over $C$. Suppose
$\mathscr E=\mathscr L\ds \sO_C^{\ds (n-1)}$ and $X=\P(\mathscr
E^*)\>C$ is the projective bundle. If $0< \deg(\mathscr
L)<\frac{2g-2}{n-1}$, then  $\P(\mathscr E^*)$ can not support
scalar-flat K\"ahler metrics, but it does admit scalar-flat
Hermitian metrics. \eproposition

\noindent As motivated by previous results, we propose the following
question.

\begin{question} Let $X$ be a compact K\"ahler manifold. Suppose $X$
has a scalar-flat Hermitian metric. Are there any geometric
conditions on $X$ which can guarantee the existence of scalar-flat
K\"ahler metrics?
\end{question}

    Finally, we classify minimal compact complex surfaces with
scalar-flat Hermitian metrics.

\btheorem\label{main3} Let $X$ be a minimal compact complex surface.
If $X$ admits a scalar-flat Hermitian metric, then $X$ must be one
of the following \bd
\item  an Enriques surface;\item a
bi-elliptic surface;\item  a K$3$ surface;
\item a $2$-torus;
\item a {Kodaira surface};
\item a ruled surface $X$ over a  curve $C$ of genus $g\geq 2$ and
$m(X)>2-2g$;
\item  a class $\mathrm{VII}_0$ surface with $b_2>0$.
 \ed\etheorem
\bremark It is  proved that surfaces in $(1)$ to $(6)$ all have
scalar-flat Hermitian metrics. On the other hand, since class
$\mathrm{VII}_0$ surfaces with $b_2>0$ are not completely
classified, we do not prove each class $\mathrm{VII}_0$ surface with
$b_2>0$ can support scalar-flat Hermitian metrics. Non-minimal
surfaces with scalar-flat Hermitian metrics will also be studied in
the sequel. \eremark

  The rest of the paper is organized as follows. In Section \ref{character}, we
 give a characterization of compact complex manifolds with
 scalar-flat Hermitian metrics and prove Theorem \ref{main1}. In
 Section \ref{ruled}, we classify ruled surfaces with scalar-flat
 Hermitian metrics and establish Theorem \ref{main2}. In Section
 \ref{minimal}, we classify minimal complex surfaces with
 scalar-flat Hermitian metrics and obtain Theorem \ref{main3}. In
 Section \ref{exs}, we give some precise examples with scalar-flat
 Hermitian metrics (Proposition \ref{main01}).\\

\textbf{Acknowledgements.} The first author would like to  thank his
advisor Professor Jian Zhou for his guidance. The second author is
very grateful to Professor
 K.-F. Liu and Professor S.-T. Yau for their support, encouragement and stimulating
discussions over  years. We would also like to thank Professors S.
Sun, A. Futaki, G. Szekelyhidi, V. Tosatti, W.-P. Zhang and X.-Y.
Zhou for some helpful discussions.

\vskip 2\baselineskip

\section{Background materials}\label{back}

\subsection{Scalar curvature and total scalar curvature on complex manifolds}
Let $(\mathscr E,h)$ be a Hermitian holomorphic vector bundle over a
complex manifold $X$ with Chern connection $\nabla$. Let
$\{z^i\}_{i=1}^n$ be the  local holomorphic coordinates
  on $X$ and  $\{e_\alpha\}_{\alpha=1}^r$ be a local frame
 of $\mathscr E$. The curvature tensor $R^{\mathscr E}\in \Gamma(X,\Lambda^{1,1}T^*_X\ts \mathrm{End}(\mathscr E))$ has components \beq R^\mathscr E_{i\bar j\alpha\bar\beta}= -\frac{\p^2
h_{\alpha\bar \beta}}{\p z^i\p\bar z^j}+h^{\gamma\bar
\delta}\frac{\p h_{\alpha \bar \delta}}{\p z^i}\frac{\p
h_{\gamma\bar\beta}}{\p \bar z^j}.\label{cu}\eeq (Here and
henceforth we sometimes adopt the Einstein convention for
summation.) If $(X,\omega_g)$ is a  Hermitian manifold, then
$(T_X,g)$ has Chern curvature components \beq R_{i\bar j k\bar
\ell}=-\frac{\p^2g_{k\bar \ell}}{\p z^i\p\bar z^j}+g^{p\bar
q}\frac{\p g_{k\bar q}}{\p z^i}\frac{\p g_{p\bar \ell}}{\p\bar
z^j}.\eeq The Chern-Ricci curvature $\mathrm{Ric}(\omega_g)$ of
$(X,\omega_g)$ is represented by $R_{i\bar j}=g^{k\bar\ell} R_{i\bar
j k\bar\ell}$.  The \emph{(Chern) scalar curvature} $s$ of
$(X,\omega_g)$ is given by \beq
s=\text{tr}_{\omega_g}\text{Ric}(\omega_g)=g^{i\bar j} R_{i\bar j}.
\eeq The \emph{total (Chern) scalar curvature} of $\omega_g$ is \beq
\int_X s\omega_g^n=n \int\text{Ric}(\omega_g)\wedge
\omega_g^{n-1},\eeq where $n$ is the complex dimension of $X$. \bd
\item A Hermitian metric $\omega_g$ is called a Gauduchon metric if
$\p\bp\omega_g^{n-1}=0$. It is proved by Gauduchon (\cite{Gau77})
that, in the conformal class of each Hermitian metric, there exists
a unique Gauduchon metric (up to constant scaling).
%

\item  A projective manifold $X$ is called uniruled if it is
covered by rational curves.
 \ed

\subsection{Positivity of line bundles.}

Let $(X,\omega_g)$ be a compact Hermitian manifold, and $\mathscr
L\>X$ be a holomorphic line bundle. \bd\item $\mathscr L$ is said to
be \emph{positive} (resp. \emph{semi-positive}) if there exists a
smooth Hermitian metric $h$ on $\mathscr L$ such that the curvature
form $R^{\mathscr L}=-\sq\p\bp\log h$ is a positive (resp.
semi-positive) $(1,1)$-form.

\item $\mathscr L$ is said to be \emph{nef}, if for any  $\eps>0$, there exists a
smooth Hermitian metric $h_\eps$ on $\mathscr L$ such that $
-\sq\p\bp\log h_\eps\geq -\eps \omega_g.$

\item $\mathscr L$ is said to be \emph{pseudo-effective}, if there exists a
(possibly) singular Hermitian metric $h$ on $\mathscr  L$ such that
 $-\sq\p\bp\log h\geq 0$
in the sense of distributions. (See \cite{Dem} for more details.)

\item $\mathscr L$ is said to be \emph{ $\Q$-effective}, if there exists some
positive integer $m$ such that $H^0(X,\mathscr  L^{\ts m})\neq 0$.

\item $\mathscr L$ is called \emph{unitary flat} if there exists a smooth Hermitian metric $h$ on $\mathscr L$ such that the curvature of
$(\mathscr L,h)$ is zero, i.e. $-\sq\p\bp\log h=0$.

 \item  The Kodaira dimension $\kappa(\mathscr  L)$ of $\mathscr  L$ is defined to
be
$$\kappa(\mathscr  L):=\limsup_{m\>+\infty} \frac{\log \dim_\C
H^0(X,\mathscr  L^{\ts m})}{\log m}$$ and the \emph{Kodaira
dimension} $\kappa(X)$ of $X$ is defined as $
\kappa(X):=\kappa(K_X)$ where the logarithm of zero is defined to be
$-\infty$. \ed
\subsection{Positivity of vector bundles}
 The points of the projective bundle $\P(\mathscr E^*)$ of
$\mathscr E\>X$ can be identified with the hyperplanes of $\mathscr
E$. Note that every hyperplane $\mathscr V$ in $\mathscr E_z$
corresponds bijectively to the line of linear forms in $\mathscr
E^*_z$ which vanish on $\mathscr V$.
 Let $\pi : \P(\mathscr E^*) \> X$ be
the natural projection. There is a tautological hyperplane subbundle
$\mathscr  S$ of $\pi^*\mathscr E$ such that $\mathscr S_{[\xi]} =
\xi^{-1}(0) \subset \mathscr E_z$ for all $\xi\in \mathscr
E_z^*\setminus\{0\}$. The quotient line bundle $\pi^*\mathscr
E/\mathscr  S$ is denoted $\sO_\mathscr E(1)$ and is called the
\emph{tautological line bundle} associated to $\mathscr E\>X$. Hence
there is an exact sequence of vector bundles over $\P(\mathscr
E^*)$, $ 0 \>\mathscr S\>\pi^*\mathscr E\>\sO_\mathscr E(1)\> 0.$ A
holomorphic vector bundle $\mathscr E\>X$ is called \emph{ample}
(resp. \emph{nef}) if the line bundle $\sO_\mathscr E(1)$ is ample
(resp. nef) over $\P(\mathscr E^*)$. (\textbf{Caution}: In general,
$\P(\mathscr E)$ and $\P(\mathscr E^*)$ are not
 isomorphic! $\sO_{\mathscr E}(1)$ is the tautological line bundle
 of $\P(\mathscr E^*)$, and $\sO_{\mathscr E^*}(1)$ is the tautological line bundle
 of $\P(\mathscr E)$.)
 A Hermitian holomorphic vector bundle $(\mathscr E,h)$
over a complex manifold $X$ is called \emph{Griffiths positive} if
at each point $q\in X$ and for any nonzero vector $v\in \mathscr
E_q$, and any nonzero vector $u\in T_qX$,  $ R^\mathscr E(u,\bar
u,v,\bar v)>0.$

\subsection{RC-positive line bundles}
Let's recall that \bdefinition A line bundle $\mathscr  L$ is called
\emph{RC-positive}  if it has a smooth Hermitian metric $h$ such
that its
 curvature $R^{(\mathscr L, h)}=-\sq\p\bp\log h$ has at least one positive eigenvalue
 everywhere.
 \edefinition

 \noindent
  In \cite[Theorem~1.4]{Yang17D}, we obtained an
 equivalent characterization for  RC-positive line
 bundles.
 \btheorem\label{equi} Let $\mathscr L$ be a holomorphic line bundle over a compact complex manifold $X$. The
 following statements  are equivalent.

 \bd \item $\mathscr L$ is  RC-positive;
 \item the dual line bundle $\mathscr L^{*}$ is not pseudo-effective.
 \ed
 \etheorem

\noindent Hence, we obtain   \bcorollary\label{flat} A line bundle
$\mathscr L$ is unitary flat if and only if neither $\mathscr L$ nor
$\mathscr L^{*}$ is RC-positive. \ecorollary \bproof It is easy to
see that $\mathscr L$ is unitary flat if and only if both $\mathscr
L$ and $\mathscr L^*$ are pseudo-effective (e.g.
\cite[Theorem~3.4]{Yang17a}). Hence, Corollary \ref{flat} follows
from Theorem \ref{equi}. \eproof

\noindent By using Theorem \ref{equi}, the classical result of
\cite[Theorem]{BDPP13} and Yau's theorem \cite{Yau78}, we obtain in
\cite[Corollary~1.9]{Yang17D} that

\btheorem\label{uniruled} A projective manifold $X$ is uniruled if
and only if $K_X^{-1}$ is RC-positive, i.e. $X$ has a smooth
Hermitian metric $\omega$ such that the Ricci curvature
$\mathrm{Ric}(\omega)$ has at least one positive eigenvalue
everywhere.

\etheorem


\vskip 2\baselineskip

\section{Characterizations of complex manifolds with scalar-flat  metrics
}\label{character}

In this section, we shall prove Theorem \ref{main1}. Let $\omega$ be
a smooth Hermitian metric on a compact complex manifold $X$. For
simplicity, we denote by $\mathscr F(\omega)$ the total (Chern)
scalar curvature of $\omega$, i.e.
$$\mathscr F(\omega)=\int_X s \omega^n=n
\int_X\text{Ric}(\omega)\wedge \omega^{n-1}.$$ Note that, when $X$
is not K\"ahler, the total  scalar curvature differs from the total
scalar curvature of the Levi-Civita connection of the underlying
Riemannian metric (e.g. \cite{LY14}). Let $\mathscr W$ be the space
of smooth Gauduchon metrics on $X$. We obtained in
\cite[Theorem~1.1]{Yang17a} a complete characterization on the image
of the total scalar curvature function $\mathscr F:\mathscr W\>\R$
following  \cite{Gau77, Mi82, La99} (see also some special cases in
\cite{Tel06, Gau77, HW12}). By Theorem \ref{equi}, we obtain the
following result.

\btheorem\label{table} The image of the total scalar function
$\mathscr F:\mathscr W\>\R$ has exactly four different cases: \bd
\item $\mathscr F(\mathscr W)=\R$ if and only if both $K_X$ and
$K_X^{-1}$ are RC-positive;
\item $\mathscr F(\mathscr W)=\R^{>0}$ if and only if $K^{-1}_X$ is RC-positive but $K_X$ is not RC-positive;
\item $\mathscr F(\mathscr W)=\R^{<0}$ if and only if $K_X$ is RC-positive  but $K^{-1}_X$ is not
RC-positive;
\item $\mathscr F(\mathscr W)=\{0\}$ if and only if $X$ is Ricci-flat; or equivalently, neither $K_X$ nor $K_X^{-1}$ is RC-positive.
\ed \etheorem \bproof We obtained in \cite[Theorem~1.1]{Yang17a}
that the image of the total scalar function $\mathscr F:\mathscr
W\>\R$ has exactly four different cases: \bd
\item $\mathscr F(\mathscr W)=\R$, if and only if neither $K_X$ nor
$K_X^{-1}$ is pseudo-effective;
\item $\mathscr F(\mathscr W)=\R^{>0}$, if and only if $K_X^{-1}$ is pseudo-effective but not unitary flat;
\item $\mathscr F(\mathscr W)=\R^{<0}$, if and only if $K_X$ is pseudo-effective but not unitary flat;
\item $\mathscr F(\mathscr W)=\{0\}$, if and only if $K_X$ is  unitary flat.
\ed By \cite[Corollary~2]{TW10}, $K_X$ is unitary flat if and only
if $X$ is Ricci-flat, i.e. there exists a Hermitian metric $\omega$
on $X$ such that $\mathrm{Ric}(\omega)=0$.  Hence Theorem
\ref{table} follows from Theorem \ref{equi} and Corollary
\ref{flat}. \eproof

\bremark It is easy to see that Theorem \ref{table} also holds for
Bott-Chern classes (\cite[Theorem~3.4]{Yang17a})

\eremark

As an application of Theorem \ref{table}, we establish Theorem
\ref{main1}, that is,

\btheorem\label{main111} Let $X$ be a compact complex manifold. Then
$X$ admits a scalar-flat Hermitian metric if and only if $X$ is
Ricci-flat, or both $K_X$ and $K_X^{-1}$ are RC-positive. \etheorem

\bproof If $X$ has  a scalar-flat Hermitian metric $\omega$, in the
conformal class of $\omega$, there exists a Gauduchon metric
$\omega_f=e^f\omega$.  Then the total scalar curvature $s_f$ of the
Gauduchon metric $\omega_f$ is \beq  s_f= n\int_X
\mathrm{Ric}(\omega_f)\wedge
\omega_f^{n-1}=n\int_X\left(\mathrm{Ric}(\omega)-n\sq\p\bp
f\right)\wedge \omega_f^{n-1}.\eeq Since $\omega_f$ is Gauduchon,
i.e. $\p\bp \omega_f^{n-1}=0$, an integration by part yields \be
s_f&=&n\int_X\mathrm{Ric}(\omega)\wedge
\omega_f^{n-1}\\&=&n\int_X\mathrm{Ric}(\omega)\wedge
e^{(n-1)f}\omega^{n-1}\\&=&\int_X e^{(n-1)f}\cdot
\mathrm{tr}_\omega\mathrm{Ric}(\omega)\cdot \omega^n. \ee Since
$\omega$ has zero scalar curvature, i.e.
$\mathrm{tr}_\omega\mathrm{Ric}(\omega)=0$, we deduce that the total
scalar curvature $s_f$ of the Gauduchon metric $\omega_f$ is zero.
By Theorem \ref{table}, we conclude that either $X$ is Ricci-flat,
or both $K_X$ and $K^{-1}_X$ are RC-positive.

  On the other hand, suppose either $X$ is
Ricci-flat, or both $K_X$ and $K^{-1}_X$ are RC-positive, by Theorem
\ref{table} again, we know $X$ has a Gauduchon metric $\omega_{{G}}$
with zero total scalar curvature. By a conformal perturbation
method, it is easy to see that there exists a Hermitian metric
$\omega$ with zero scalar curvature (e.g.
\cite[Lemma~3.2]{Yang17a}). Indeed, let $s_{G}$ be the
 scalar curvature of $\omega_{{G}}$. It is well-known (e.g.\cite{Gau77}£¬
or \cite[Theorem~2.2]{CTW16}) that the following equation \beq s_G-
\text{tr}_{\omega_G}\sq\p\bp f=0\eeq has a solution $f\in
C^\infty(X)$ since $\omega_G$ is Gauduchon and its total scalar
curvature $\int_Xs_G\omega^n_G$ is zero. Let $ \omega=
e^{\frac{f}{n}}\omega_G$. Then the scalar curvature $s$ of $\omega$
is, \be s&=&\text{tr}_{\omega }\text{Ric}(\omega)=-\text{tr}_{\omega
}\sq \p\bp\log(
\omega^n)\\&=&-e^{-\frac{f}{n}}\text{tr}_{ \omega_G}\sq\p\bp\log (e^f \omega_G^n) \\
&=&-e^{-\frac{f}{n}}\left(s_G-\text{tr}_{\omega_G}\sq\p\bp
f\right)\\&=&0.\ee The proof of Theorem \ref{main1} is completed.
\eproof

\emph{The proof of Corollary \ref{main0}.} It is a special case of
Theorem \ref{main1} since K\"ahler manifolds with unitary flat $K_X$
are  K\"ahler Calabi-Yau.\qed

\bcorollary Let $X$ be a compact K\"ahler manifold. Suppose $X$ has
a scalar-flat Hermitian metric, or a Gauduchon metric with zero
total scalar curvature. If $K_X$ or $K_X^{-1}$ is pseudo-effective,
then $X$ is a K\"ahler Calabi-Yau manifold. \ecorollary

\vskip 2\baselineskip

\section{Projective bundles with scalar-flat
metrics}\label{technique}
In this section, we prove the following result.

 \btheorem\label{key}
Let $\mathscr E$ be a nef vector bundle of rank $r\geq 2$ over a
smooth curve $C$ with genus $g\geq 2$ and $X=\P(\mathscr E)$.  If
$0\leq \deg(\mathscr E)<2g-2$, then both  $K_X$ and  $K_X^{-1}$ are
RC-positive. In particular, $X$ has scalar-flat Hermitian metrics.
\etheorem

Let's recall some elementary settings. Suppose $\dim_\C Y=n$ and
$r=\mathrm{rank}(\mathscr E)$. Let $\pi$ be the projection
$\P(\mathscr E^*)\> Y$ and $\mathscr L=\sO_\mathscr E(1)$. Let
$(e_1,\cdots, e_r)$ be the local holomorphic frame on $\mathscr E$
and the dual frame on $\mathscr E^*$ is denoted by $(e^1,\cdots,
e^r)$. The corresponding holomorphic coordinates on $\mathscr E^*$
are denoted by $(W_1,\cdots, W_r)$.  If $(h_{\alpha\bar\beta})$ is
the matrix representation of a smooth  metric $h^{\mathscr E}$ on
$\mathscr E$ with respect to the basis $\{e_\alpha\}_{\alpha=1}^r$,
then the induced Hermitian metric $h^\mathscr L$ on $\mathscr L$ can
be written as $ h^\mathscr L=\frac{1}{\sum
h^{\alpha\bar\beta}W_\alpha\bar W_\beta}$.
 The curvature of $(\mathscr L,h^\mathscr L)$ is \beq
R^{\mathscr L}=\sq\p\bp\log\left(\sum
h^{\alpha\bar\beta}W_\alpha\bar W_\beta\right)
\label{inducedcurvature}\eeq where $\p$ and $\bp$ are operators on
the total space $\P(\mathscr E^*)$. We fix a point $p\in \P(\mathscr
E^*)$, then there exist local holomorphic coordinates
 $(z^1,\cdots, z^n)$ centered at point $q=\pi(p)\in Y$ and local holomorphic basis $\{e_1,\cdots, e_r\}$ of $\mathscr E$ around $q$ such that
 \beq h_{\alpha\bar\beta}=\delta_{\alpha\bar\beta}-R^{\mathscr E}_{i\bar j \alpha\bar\beta}z^i\bar z^j+O(|z|^3) \label{normal}\eeq
Without loss of generality, we assume $p$  is the point $(0,\cdots,
0,[a_1,\cdots, a_r])$ with $a_r=1$. On the chart $U=\{W_r=1\}$ of
the fiber $\P^{r-1}$, we set $w^A=W_A$ for $A=1,\cdots, r-1$. By
formula (\ref{inducedcurvature}) and (\ref{normal}) \beq R^{\mathscr
L}(p)=\sq\sum R^{\mathscr E}_{i\bar j\alpha\bar\beta}\frac{a_\beta
\bar a_\alpha}{|a|^2}dz^i\wedge d\bar z^j+\omega_{\mathrm{FS}}
\label{induced} \eeq where
$|a|^2=\sum\limits_{\alpha=1}^r|a_\alpha|^2$ and
$\omega_{\mathrm{FS}}=\sq
\sum\limits_{A,B=1}^{r-1}\left(\frac{\delta_{AB}}{|a|^2}-\frac{a_B\bar
a_A}{|a|^4}\right)dw^A\wedge d\bar w^B$ is the Fubini-Study metric
on the fiber $\P^{r-1}$.

\blemma\label{RCample}  If $\mathscr E$ is Griffiths-positive, then
$\sO_{\mathscr E^*}(-1)$ is RC-positive. \elemma \bproof It follows
from formula (\ref{induced}). Indeed, by (\ref{induced}), the
induced metric on $\sO_{\mathscr E^*}(-1)$ over $\P(\mathscr E)$ has
curvature form $$ R^{\sO_{\mathscr E^*}(-1)}=-\left(\sq\sum
R^{\mathscr E^*}_{i\bar j\alpha\bar\beta}\frac{a_\beta \bar
a_\alpha}{|a|^2}dz^i\wedge d\bar z^j+\omega_{\mathrm{FS}}\right). $$
On the other hand, $R^{\mathscr E^*}=-\left(R^{\mathscr E}\right)^t$
and so $$ R^{\sO_{\mathscr E^*}(-1)}=\sq\sum R^{\mathscr E}_{i\bar
j\alpha\bar\beta}\frac{a_\alpha \bar a_\beta}{|a|^2}dz^i\wedge d\bar
z^j-\omega_{\mathrm{FS}}. $$ Hence,
 $\sO_{\mathscr
E^*}(-1)$ is RC-positive  if $(\mathscr E, h^{\mathscr E})$ is
Griffiths-positive. \eproof

\blemma\label{RCnef} If $\mathscr E$ is a nef vector bundle over a
smooth curve $C$. Then for any ample line bundle $\mathscr A$ over
$C$ and any $k\geq 0$, $\sO_{\mathscr E^*}(-k)\ts \pi^*(\mathscr A)$
is RC-positive.

 \bproof It is easy to see that $\mathrm{Sym}^{\ts
k}\mathscr E \ts \mathscr A$ is an ample vector bundle over $C$. By
\cite{CF90}, $\mathrm{Sym}^{\ts k}\mathscr E\ts \mathscr A$ has a
smooth Griffiths-positive metric. In particular, by Lemma
\ref{RCample}, the dual tautological line bundle \beq
\sO_{\mathrm{Sym}^{\ts k}\mathscr E^*\ts \mathscr A^*}(-1)\eeq is
RC-positive. More precisely, the base curve $C$ direction is a
positive direction of the curvature tensor of
$\sO_{\mathrm{Sym}^{\ts k}\mathscr E^*\ts \mathscr A^*}(-1)$. On the
other hand, we have the following commutative diagram
$$\CD
  \sO_{\mathscr E^*}(-k)\ts \pi^*(\mathscr A) @>>> \sO_{\mathrm{Sym}^{\ts k}\mathscr E^*}(-1)\ts \pi_k^*(\mathscr A) @>>>\sO_{\mathrm{Sym}^{\ts k}\mathscr E^*\ts \mathscr
A^*}(-1)\\
  @V  VV @V  VV  @V  VV \\
  \P(\mathscr E) @>\nu_k>> \P(\mathrm{Sym}^{\ts k}\mathscr E) @>i>> \P(\mathrm{Sym}^{\ts k}\mathscr E\ts \mathscr A)\\
  @V_{\pi}  VV @V_{\pi_k}  VV @V  VV \\
  C @>f>> C@>f>>C,
\endCD
$$
where $\nu_k:\mathscr E\>\mathrm{Sym^{\ts k}}\mathscr E$ is the
$k$-th Veronese map,  $f=\mathrm{Identity}$ and $i$ is an
isomorphism. It is easy to see that $\sO_{\mathscr E^*}(-k)\ts
\pi^*(\mathscr A)$ is RC-positive, i.e., the induced curvature has a
positive direction along the base $C$ direction. \eproof \elemma

\noindent\emph{The proof of Theorem \ref{key}.} By using the
projection formula on $X=\P(\mathscr E)$, $$ K_X=\sO_{\mathscr
E^*}(-n)\ts \pi^*(K_C\ts \det \mathscr E^*), $$ where $\pi:X\>C$ is
the projection. If $\deg(\mathscr E)<2g-2=\deg(K_C)$, then $\deg
(K_C\ts \det \mathscr E^*)>0$ and so $K_C\ts \det \mathscr E^*$ is
ample. By Lemma \ref{RCnef}, $K_X$ is RC-positive. On the other
hand, by Theorem \ref{uniruled}, it is easy to see that $K_X^{-1}$
is RC-positive. Hence, by Theorem \ref{main1}, $X$ has scalar-flat
Hermitian metrics. \qed

\vskip 2\baselineskip

\section{Classification of ruled surfaces with scalar-flat Hermitian
metrics}\label{ruled}

In this section, we classify ruled surfaces with scalar-flat
Hermitian metrics and prove Theorem \ref{main2}. It is well-known
that any ruled surface $X$ can be written as a projective bundle
$\P(\mathscr E)$ where $\mathscr E$ is a rank two vector bundle over
a smooth curve $C$ with genus $g$. Moreover, two ruled surfaces
$\P(\mathscr E)$ and $\P(\mathscr E')$ are isomorphic if and only if
$\mathscr E\cong \mathscr E'\otimes \mathscr L$ for some line bundle
$\mathscr L$ over $C$.  Since  $\mathscr E$ has rank two and $X\cong
\P(\mathscr E)\cong \P(\mathscr E^*)$,  we shall use projection
formulas $$K_X=\sO_{\mathscr E}(-2)\ts \pi^*(K_C\ts \det \mathscr
E),\ \ \pi:\P(\mathscr E^*)\>C$$ and  $$K_X=\sO_{\mathscr
E^*}(-2)\ts \pi^*(K_C\ts \det \mathscr E^*),\ \ \pi:\P(\mathscr
E)\>C$$ alternatively.

 When $g=0$, $C\cong \P^1$ and each rank two vector bundle can be written as $\mathscr E=\sO_{\P^1}(a)\ds
 \sO_{\P^1}(b)$. We can write a ruled surface over $\P^1$ as $X=\P(\sO_{\P^1}(-k)\ds
\sO_{\P^1})$.

\bproposition\label{rational} Let $X=\P(\sO_{\P^1}(-k)\ds
\sO_{\P^1})$  be a Hirzebruch surface. Then the anti-canonical line
bundle $K_X^{-1}$ is effective and $X$ has no scalar-flat Hermitian
metrics. \eproposition \bproof Let $\mathscr E=\sO_{\P^1}(k)\ds
\sO_{\P^1}$ and $X=\P(\mathscr E^*)$. We have $K^{-1}_X=\sO_\mathscr
E(2)\ts \pi^*(\sO_{\P^1}(2-k)).$ By the direct image formula (e.g.
\cite[p.90]{Laza1}), we have \be H^0(X,K_X^{-1})&=&
H^0(X,\sO_\mathscr E(2)\ts
\pi^*(\sO_{\P^1}(2-k))\\
&=&H^0(\P^1, \text{Sym}^{\ts 2}\mathscr E\ts \sO_{\P^1}(2-k))\\
&=&H^0(\P^1,\sO_{\P^1}(k+2)\ds
\sO_{\P^1}(2)\ds\sO_{\P^1}(2-k))\\
&\neq& 0\ee for any $k$. Therefore, $K_X^{-1}$ is effective and
$K_X$ is not RC-positive. By Theorem \ref{main1}, $X$ has no
scalar-flat Hermitian metrics. \eproof

\btheorem\label{elliptic} Let $X=\P(\mathscr E^*)\>C$  be a
projective bundle over an elliptic curve $C$ where $\mathscr E\>C$
is a rank two vector bundle. Then the $K_X$ is not RC-positive and
$X$ has no scalar-flat Hermitian metrics. \etheorem \bproof We
divide the proof into three different
cases.\\

 \noindent \emph{Case $1$.} Suppose $\mathscr E $ is
indecomposable and $\deg \mathscr E=0$. A well-known result of
Atiyah asserts that an indecomposable vector bundle over an
 elliptic curve is semi-stable and so $\mathscr E$ is semi-stable (e.g.
 \cite[Appendix~A]{Tu93}). On the other hand, a semi-stable vector bundle
over a curve is nef
 if  $\deg(\mathscr E)\geq 0$ (e.g. \cite[Theorem~6.4.15]{Laza1}).
 Hence $\mathscr E$ is nef. By using the projection formula,
 \beq
K^{-1}_X=\sO_{\mathscr E}(2)\ts \pi^*(K^{-1}_C\ts \det \mathscr
E^*)=\sO_{\mathscr E}(2)\ts \pi^*(\det \mathscr E^*) \label{pp}\eeq
we deduce $K_X^{-1}$ is nef.\\
%
%
%
%
%

\noindent \emph{Case $2$}. Suppose $\mathscr E$ is indecomposable
and $\deg(\mathscr E)\neq 0$. There exists an \'etale base change
$f:C'\>C$ of degree $k$ where $k$ is an integer such that $2|k$, and
$C'$ is also an elliptic curve. Suppose $X'=\P(f^*\mathscr E^*)$,
then we have the commutative diagram
$$\xymatrix{
  X' \ar[d]_{\pi'} \ar[r]^{f'}
                & X\ar[d]^{\pi}  \\
  C'  \ar[r]^{f}
                & C.            }$$
                 Let $\ell$ be an integer defined as \beq
\ell=\frac{\deg(f^*\mathscr E)}{2}=\frac{k\deg(\mathscr E)}{2},\eeq
 and $\mathscr F$ be a line bundle over $Y$ such that $\deg(\mathscr F)=-\ell$.
   Now we set $$\tilde {\mathscr E}=f^*\mathscr E\ts \mathscr F,$$ then $\deg(\tilde {\mathscr E})=0$.
 Since $\mathscr E$ is indecomposable, it is semi-stable. Therefore $f^*\mathscr E$ is
 semi-stable (e.g. \cite[Lemma~6.4.12]{Laza1}) and so $\tilde{\mathscr
 E}
 $ is semi-stable. Therefore, $\tilde{\mathscr E}$ is nef since $\deg(\tilde{\mathscr E})=0$. By
 projection formula again, we have
 $$K^{-1}_{X'}=\sO_{\tilde {\mathscr E}}(2)\ts \pi^*(\det \tilde {\mathscr E}).$$
We deduce $K^{-1}_{X'}$ is  nef. Hence $K_X^{-1}$ is nef.\\

 \noindent \emph{Case $3$}. If $\mathscr E$ is decomposable, then
 there exits a  line bundle $\mathscr L$ such that $$\mathscr E=\mathscr L\ds (\mathscr L^{-1}\ts \det \mathscr E).$$
By the projection formula (\ref{pp}) again, we have \be
H^0(X,K_X^{-1})&=& H^0(X,\sO_\mathscr E(2)\ts \pi^*(\det\mathscr
E^*))\cong
H^0(C, \text{Sym}^{\ts 2}\mathscr E\ts \det\mathscr E^*))\\
&=&H^0\left(C,\left(\mathscr L^2\ts \det\mathscr E^*\right)\ds
\sO_{C}\ds\left(\mathscr L^{-2}\ts \det\mathscr E\right)\right)\\
&\neq& 0\ee So $K^{-1}_X$ is effective.\\

In summary, we conclude that the anti-canonical line bundle
$K_X^{-1}$ is pseudo-effective, i.e. $K_X$ is not RC-positive. By
Theorem \ref{main1}, $X$ has no scalar-flat Hermitian metrics.
 \eproof

Finally, we deal with ruled surfaces over curves of genus $g\geq 2$.
For a rank two vector bundle $\mathscr E$ over a curve $C$, in
general, it is not clear whether $\mathscr E$ has an extension by
$\sO_C$:
\begin{equation}
0\rw \sO_C\rw \mathscr E\rw \ \mathscr F\rw0
\end{equation}
where $\mathscr F$ is a coherent sheaf over $C$. However, one can
obtain such an extension for $\mathscr E\otimes \mathscr L$ where
$\mathscr L$ is some suitable line bundle. This enables us to make
the following definition (see \cite[p.121-p.124]{Fri98} for more
details).

 \bdefinition Let $\mathscr E$ be a rank two vector bundle over a smooth curve $C$.  The number $m(\mathscr E)$
is defined to be the minimal degree of $\mathscr E\otimes \mathscr
L$ where there exists a sheaf extension of $\mathscr E\otimes
\mathscr L$: \beq 0\rw \sO_C\rw \mathscr E\otimes \mathscr L\rw
\mathscr F\rw0 \label{extension}\eeq      for some line bundle
$\mathscr L$ over $C$. \edefinition

\noindent It is easy to see that for a sufficiently ample line
bundle $\mathscr L$,  $H^0(C, \mathscr E\ts \mathscr L)\neq 0$ and a
global section of $\mathscr E\ts \mathscr L$ gives an extension
(\ref{extension}). Hence, $m(\mathscr E)$ is well-defined.  It is
obvious that $m(\mathscr E)=m(\mathscr E\ts \tilde {\mathscr L})$
for any line bundle $\tilde {\mathscr L}$. Nagata proved in
\cite[Theorem~1]{Nag70} (see also \cite[p.~123]{Fri98}) that

\btheorem\label{Nagata} $m(\mathscr E)\leq g$.\etheorem

\noindent (Note that, in \cite[p.~123]{Fri98}, the notion
$e(\mathscr E)$ is exactly $-m(\mathscr E)$.)\\

As we pointed out before, any ruled surface $X$ can be written as a
projective bundle $\P(\mathscr E)$ and two ruled surfaces
$\P(\mathscr E)$ and $\P(\mathscr E')$ are isomorphic if and only if
$\mathscr E\cong \mathscr E'\otimes \mathscr L$ for some line bundle
$\mathscr L$, then we can define $m(X)$ by $m(\mathscr E)$ for any
ruled surface $X=\P(\mathscr E)$.\\

 One can see that the definition of $m(\mathscr E)$ is related to stability of coherent sheaves. If $m(\mathscr E)>0$,
 then  $\mathscr E$ is stable. Indeed, for any rank one sub-sheaf  $\mathscr L$ of $\mathscr E$, we have the short exact sequence:
 $$0\rw \mathscr L\rw \mathscr E\rw \mathscr F\rw 0.$$
  Since $\mathscr E$ is torsion free, $\mathscr L$ is torsion free and we know $\mathscr L$ is a line
  bundle. Therefore,
 $$0 \rw \mathscr \sO_C \rw \mathscr E\otimes \mathscr L^{-1} \rw \mathscr F\otimes \mathscr L^{-1}\rw 0.$$
By the definition of $m(\mathscr E)$, we have
 $ \deg(\mathscr E \otimes \mathscr L^{-1})\geq m(\mathscr E)>0$
which is equivalent to $\deg\mathscr L<\frac{\deg \mathscr E}{2}.$
This implies $\mathscr E$ is stable. Conversely, if $\mathscr E$ is
stable, by a similar argument, we can conclude $m(\mathscr E)>0$.
Hence, we obtain a fact pointed out in \cite[Proposition~12,
p.~123]{Fri98}.

\bproposition\label{stable} If $\mathscr E$ is a rank two vector
bundle over a Riemann surface $C$, then $\mathscr E$ is stable if
and only if $m(\mathscr E)>0$. \eproposition

\noindent \emph{The proof of Theorem \ref{main2}.} Let $X$ be a
ruled surface which can support scalar-flat
 Hermitian metrics. We can write $X=\P(\mathscr E_o)$ for some rank
 $2$ vector bundle $\mathscr E_o$ over a smooth curve $C$. Note that, since $\mathscr E_o$ has
 rank $2$, $\mathscr E_o\cong \mathscr E_o^*\ts \det\mathscr E_o$ and so
 $X\cong \P(\mathscr E_o)\cong \mathscr \P(\mathscr E_o^*)$. By
 Proposition \ref{rational} and Theorem \ref{elliptic}, we know the
 genus $g(C)\geq 2$. On the other hand, by the above discussion, we can write $X=\P(\mathscr E)$ where $\deg(\mathscr
 E)=m(X)$ and $\mathscr E$ has an extension
  \beq 0\>\sO_C\>\mathscr E\>\mathscr F\>0.\label{extension1}\eeq
Hence, $\deg(\mathscr E)=\deg(\mathscr F)=m(X)$.\\

 (1).  If $m(X)=\deg \mathscr F\leq 2-2g$,  $X\cong \P(\mathscr E^*)\cong \P(\mathscr E)$ has
no scalar-flat Hermitian metrics. Indeed, we consider $X=\P(\mathscr
E^*)$.  By the exact sequence (\ref{extension1}), we have
$$0\>H^0(C,\sO_C)\>H^0(C,\mathscr E)\>\cdots. $$ Therefore,
$H^0(C,\mathscr E)\neq 0$. By the Le Potier isomorphism
(\cite{LeP75}), we have
$$H^0(\P(\mathscr E^*),\sO_{\mathscr E}(1))\cong H^0(C,\mathscr
E)\neq 0.$$ Hence, $\sO_{\mathscr E}(1)$ is effective and so it is
pseudo-effective. On the other hand, since $\deg(\mathscr E)\leq
2-2g=-\deg(K_C)$, we deduce  $K_C^{-1}\ts \det \mathscr E^*$ is
semi-positive. By the projection formula $K_X^{-1}=\sO_{\mathscr
E}(2)\ts \pi^*(K_C^{-1}\ts \det \mathscr E^*),$  we know $K_X^{-1}$
is pseudo-effective. By Theorem \ref{equi}, $K_X$ is not
RC-positive. By Theorem
\ref{main1}, $X$ has no scalar-flat Hermitian metrics.\\

(2). If $2-2g<m(X)=\deg(\mathscr E)=\deg(\mathscr F)\leq 0$, we know
$0\leq \deg(\mathscr E^*)<2g-2$. Since $\sO_C$ and $\mathscr F^*$
are nef, by the dual exact sequence of (\ref{extension1}), $$
0\>\mathscr F^* \>\mathscr E^*\>\sO_C\>0,$$ we deduce $\mathscr E^*$
is nef with $0\leq \deg(\mathscr E^*)<2g-2$. By Theorem \ref{key},
$X\cong \P(\mathscr E^*)$ can support scalar-flat Hermitian
metrics.\\

(3). If $0<m(X)=\deg(\mathscr E)=\deg(\mathscr F)<2g-2$, by the
exact sequence (\ref{extension1}), $\mathscr E$  is nef with $0<
\deg(\mathscr E)<2g-2$.  By Theorem \ref{key}, $X\cong \P(\mathscr
E)$ admits scalar-flat Hermitian metrics. Note that $\P(\mathscr
E)\cong \P(\mathscr E^*)$.\\

(4). Suppose $m(X)\geq 2g-2$. By  Theorem \ref{Nagata}, $m(X)\leq
g$. Hence, in this case, we have $g=2$ and $m(X)=\deg(\mathscr
E)=2$. We work on $X=\P(\mathscr E)$. By Proposition \ref{stable},
$\mathscr E$ is a stable vector bundle and $\deg(\mathscr E)=2$. By
(\cite[Theorem~6.4.15]{Laza1}), we know $\mathscr E$ is an ample
vector bundle over a smooth curve. According to \cite{CF90},
$\mathscr E$ has a smooth Griffiths-positive metric.  By using Lemma
\ref{RCample}, $\sO_{\mathscr E^*}(-1)$ is RC-positive. By the
projection formula again, we have
$$K_X=\sO_{\mathscr E^*}(-2)\ts \pi^*(K_C\ts \det\mathscr E^*).$$
Since $\deg(K_C)=\deg(\mathscr E)=2$, we know $K_C\ts \det\mathscr
E^*$ and $\pi^*(K_C\ts \det\mathscr E^*)$ are unitary flat. Hence,
we deduce $K_X$ is RC-positive. Since $X$ is uniruled, by Theorem
\ref{uniruled}, $K_X^{-1}$ is RC-positive. Then we can apply Theorem
\ref{main1} and assert that $X$ has  scalar-flat Hermitian
metrics.\\

 In summary, we prove that  a ruled surface $X$ over a smooth curve $C$ admits scalar-flat
 Hermitian metrics if and only if $g(C)\geq 2$ and
 $m(X)>2-2g.$ The proof of Theorem \ref{main2} is completed.
 \qed

\vskip 2\baselineskip

\section{Classification of minimal surfaces with scalar-flat Hermitian
metrics}\label{minimal}
 In this section, we classify minimal surfaces with
scalar-flat Hermitian metrics and prove Theorem \ref{main3}.

 \bproposition Let $X$ be a compact complex manifold. If $X$ admits a
 scalar-flat Hermitian metric, then the Kodaira dimension $\kappa(X)=0$ or
 $\kappa(X)=-\infty$.
 \eproposition
 \bproof According to the proof of Theorem \ref{main1}, if $X$ admits a
 scalar-flat Hermitian metric, then $X$ has a Gauduchon metric with
 zero total scalar curvature. By Theorem
 \cite[Theorem~1.4]{Yang17a},
 $\kappa(X)=0$ or
 $\kappa(X)=-\infty$.
\eproof

If $X$ is a minimal surface with Kodaira dimension $\kappa(X)=0$,
$X$ is exactly one of the following (e.g. \cite{BHPV04}) \bd\item an
Enriques surface;\item a bi-elliptic surface;\item a K$3$ surface;
\item a torus;
\item a Kodaira surface.
\ed In this case, it is well-known that $X$ has torsion canonical
line bundle, i.e. $K_X^{\ts 6}=\sO_X$ (e.g. \cite[p.~244]{BHPV04}).
Hence, $X$ admits scalar-flat Hermitian metrics.

If $X$ is a minimal surface with Kodaira dimension
$\kappa(X)=-\infty$, then $X$ lies in one of the following classes:

 \bd \item minimal rational surfaces;

 \item ruled surfaces of genus $g\geq 1$;

 \item minimal  surfaces of class $\mathrm{VII}_0$.
 \ed

 Minimal rational surfaces are either $\P^2$ or  Hirzebruch surfaces.
Hence, by Proposition \ref{rational}, they can not support
scalar-flat Hermitian metrics.

If $X$ is a minimal ruled surfaces of genus $g\geq 1$, by Theorem
\ref{main3},  $X$ has a scalar-flat Hermitian metric if and only if
$g\geq 2$ and $m(X)>2-2g$.

If $X$ is a minimal surface of class $\mathrm{VII}_0$, then $X$ is
one of the following

 \bd

\item[$\bullet$] class $\mathrm{VII}_0$ surfaces with $b_2>0$;

\item[$\bullet$] Inoue surfaces: a class $\mathrm{VII}_0$ surface has $b_2=0$ and
contains no curves;

\item[$\bullet$] Hopf surfaces: its universal covering is $\C^2-\{0\}$, or
equivalently a class $\mathrm{VII}_0$ surface has $b_2=0$ and
contains a curve.

\ed

\noindent According to the proof of \cite[Remark~4.2]{Tel06} (see
also \cite{TW13} or \cite[Theorem~5.1]{HLY18}), we know Inoue
surfaces all have $K_X$ semi-positive but not unitary flat, and so
it can not support scalar-flat Hermitian metrics. Similarly, it is
proved in \cite[Remark~4.3]{Tel06}, all Hopf surfaces have
semi-positive anti-canonical bundle, and so it has no scalar-flat
Hermitian metrics. For class $\mathrm{VII}_0$ surfaces with $b_2>0$,
they are not completely classified, and it is possible that some of
them can support scalar-flat Hermitian metrics (see the discussion
in \cite[p.~977-p.~979]{Tel06}). The proof of Theorem \ref{main3} is
completed.

\vskip 2\baselineskip

\section{Examples}\label{exs}

In this section, we exhibit several examples on ruled manifolds with
scalar-flat Hermitian metrics. As a straightforward application of
Theorem \ref{main2}, we get the following result.

 \bcorollary  Let $\mathscr L\>C$ be a line bundle over
a smooth curve of genus $g$ and $X=\P(\mathscr L\ds \sO_C)$. Then
$X$ has a scalar-flat Hermitian metric if and only if $g\geq 2$ and
$|\deg(\mathscr L)|< 2g-2$.
 \ecorollary

\noindent We can also  construct higher dimensional ruled manifolds
with scalar-flat metrics.

 \btheorem\label{ex} Let $C$ be a smooth
curve with genus $g\geq
 2$ and $\mathscr L$ be a line bundle over $C$. Suppose
$\mathscr E=\mathscr L\ds \sO_C^{\ds (n-1)}$ and $X=\P(\mathscr
E^*)\>C$ is the projective bundle.  If $0\leq \deg(\mathscr
L)<\frac{2g-2}{n-1}$, then both $K_X$ and $K_X^{-1}$ are
RC-positive.

\etheorem

\bproof By using the projection formula, we know \beq
K_X=\sO_{\mathscr E}(-n)\ts \pi^*(K_C\ts \det \mathscr E),
\label{p}\eeq where $\pi:X\>C$ is the projection.
   Fix an arbitrary smooth
Hermitian metric $h^{\mathscr L}$ on $\mathscr L$ and the trivial
metric on $\sO_C$. Let $\{z\}$ be the local holomorphic coordinate
on $C$. The curvature form of $(\mathscr L,h^{\mathscr L})$ is \beq
R^{\mathscr L}=-\sq\p\bp\log h^{\mathscr L}=\sq \kappa dz\wedge
d\bar z.\eeq Similarly, fix a smooth metric $h^{K_C}$ on  $K_C$, and
its curvature form is \beq R^{K_C}=-\sq\p\bp\log h^{K_C}= \sq \gamma
dz\wedge d\bar z.\eeq Hence, $\mathscr E$ has the curvature form
\beq R^{\mathscr E}= \sq \kappa dz\wedge d\bar z\ts e^1\ts
e^1+\sum_{i=2}^n \sq \cdot 0\cdot dz\wedge d\bar z\ts e^i\ts e^i,
\eeq where $e^1=e_\mathscr L$ is the local frame of $\mathscr L$ and
for $i\geq 2$, $e^i=e$  is the local holomorphic frame on $\sO_C$
with the order in the direct sum  $\mathscr E=\mathscr
L\ds\sO_C^{\ds(n-1)}$. Therefore, by (\ref{induced}), $\sO_{\mathscr
E}(-n)$ has the curvature form at some point $$ R^{\sO_{\mathscr
E}(-n)}= \sq\left(-n\kappa \frac{|a_1|^2}{|a|^2}dz\wedge d\bar
z\right)-n\omega_{\mathrm{FS}}.
$$ Hence, by formula (\ref{p}), the curvature of $K_X$ is given by
$$ R^{K_X}= \sq\left(\left((\kappa+\gamma)-n\kappa
\frac{|a_1|^2}{|a|^2}\right)dz\wedge d\bar
z\right)-n\omega_{\mathrm{FS}}.
 $$ Since $\deg(\mathscr L)\geq 0$, we can choose the smooth
metric $h^{\mathscr L}$ such that its curvature is semi-positive,
i.e. $\kappa\geq 0$. Therefore, \beq  R^{K_X}\geq
\sq\left(\left(\gamma-(n-1)\kappa \right)dz\wedge d\bar
z\right)-n\omega_{\mathrm{FS}}. \label{RC} \eeq
  The condition $0\leq
\deg(\mathscr L)<\frac{2g-2}{n-1}$ implies $\deg\left(K_C\ts
\mathscr L^{1-n}\right)>0$. Therefore, we can choose the Hermitian
metric $h^{K_C}$ on $K_C$ such that $h^{K_C}\ts (h^{\mathscr
L})^{1-n}$ has positive curvature, i.e. $$\gamma-(n-1)\kappa>0.$$ By
(\ref{RC}), we know the curvature of $K_X$ is positive along the
base direction, i.e., $K_X$ is RC-positive. The RC-positivity of
$K_X^{-1}$ follows from Theorem \ref{uniruled}.   \eproof

\bexample Let $n\geq 2$ be an integer. Let $C$ be a smooth curve of
degree $d\geq n+3$ in $\P^2$. It is easy to see that
$\mathrm{deg}(\sO_C(1))=d$ and $C$ is a curve of genus \beq
g=\frac{(d-1)(d-2)}{2}. \eeq
 Let $\mathscr L=\sO_C(1)$ and  $\mathscr
E=\mathscr L\ds\sO_C^{\ds(n-1)}$ and $X:=\P(\mathscr E^*)\>C$ be the
projective bundle. Note that $\dim_\C X=n$. Then
$$\frac{2g-2}{n-1}=\frac{d(d-3)}{n-1}\geq \frac{d\cdot n}{n-1}>d=\deg(\mathscr L)>0.$$
Hence, the pair $(X, C, \mathscr L, \mathscr E)$ satisfies the
conditions in Theorem \ref{ex}. In particular, both $K_X$ and
$K_X^{-1}$ are RC-positive. \eexample

\emph{The proof of Proposition \ref{main01}.} By Theorem \ref{ex}
and Theorem \ref{main1}, $X$ admits a scalar-flat Hermitian metric.
On the other hand, by \cite[Theorem~1]{ACGT11}, $X$ has no
scalar-flat K\"ahler metrics since $\mathscr E=\mathscr
L\ds\sO_C^{\ds(n-1)}$ is not polystable. \qed

\end{document}